\documentclass{elsarticle}
\newcommand{\erd}{Erd\H{o}s}
\let\ts=\thinspace

\usepackage{hyperref,amsmath}
\usepackage{isabelle,isabellesym}
\usepackage{amssymb}
\usepackage[T1]{fontenc}
\usepackage[utf8]{inputenc}             

\journal{Annals of Pure and Applied Logic}

\newcommand\narrows{
  \mathrel{\mkern2.7mu\not\mkern-2.7mu\longrightarrow}}
\DeclareMathOperator{\tp}{tp}

\begin{document}

\begin{frontmatter}

\title{A Formalised Theorem in the Partition Calculus}

\author{Lawrence C. Paulson FRS}
\address{Computer Laboratory, University of Cambridge\\
   15 JJ Thomson Avenue, Cambridge CB3 0FD, UK}


\begin{abstract}
A paper on ordinal partitions by \erd{} and Milner~\cite{erdos-theorem-partition} has been formalised using the proof assistant Isabelle/HOL, augmented with a library for Zermelo--Fraenkel set theory.
The work is part of a project on formalising the partition calculus. The chosen material is particularly appropriate in view of the substantial corrections~\cite{erdos-theorem-partition-corr} later published by its authors, illustrating the potential value of formal verification.
\end{abstract}


\end{frontmatter}
 
\section{Introduction}

Formal logic was developed to strengthen the foundations of mathematics.
Whitehead and Russell's magnum opus~\cite{principia} may have been intended to show that all of mathematics could be formalised,
but rather suggested the opposite: they managed to prove
$1+1 = 2$ only on page 360.
The highly formal mathematics of Bourbaki has been sharply criticised by A.~R.~D. Mathias, who (among other things) points out \cite{Mathias2002} that 
their definition of the number 1 expands to some $4.5\times 10^{12}$ symbols.
While many mathematicians are indifferent to logic, 
Mathias' criticism is particularly trenchant in that he himself is a logician.
For all that, researchers today are trying to formalise mathematics using formal deductive logic, 
with the help of software called \textit{interactive theorem provers} (or \textit{proof assistants}).

The field of automated theorem proving was initiated by logicians and philosophers such as Martin Davis and Hilary Putnam (whose early work~\cite{davis-putnam} eventually led to today's powerful DPLL procedure for propositional logic) and Alan Robinson (the father of the resolution method for theorem proving in first-order logic~\cite{robinson65}). One might have expected the next step to be the automation of set theory,
but instead the field took a sharp turn: away from full automation to interaction, with a  focus on problems in computer science.
A milestone was Michael J. C. Gordon's focus on hardware verification and his choice of higher-order logic~\cite{mgordon86}, a choice that logicians would not have made, as it had ``no coherent well-established theory''~\cite[p.\ts241]{van-benthem-higher} compared with first-order logic. Soon, researchers around the world were experimenting with Gordon's interactive theorem prover, HOL~\cite{mgordon-hol}. 
Other implementations of higher-order logic soon appeared, such as HOL Light \cite{harrison-hol-light} and Isabelle/HOL~\cite{isa-tutorial}.

The scope of higher-order logic turned out to be much greater than hardware verification. Many researchers
turned to the formalisation of well-known mathematical results~\cite{harrison-pnt}
and even to the verification of a contested result: Hales' computer-assisted proof of the Kepler conjecture~\cite{hales-formal-Kepler}.
A separate strand of research based upon constructive type theories also led to the
formalisation of deep mathematical results, such as the odd order theorem~\cite{gonthier-oot}.
The importance of such achievements has been recognised by the Isaac Newton Institute's
programme entitled \textit{Computer-aided Mathematical Proof} (2017)
and in the 2020 Mathematics Subject Classification \cite{dunne-mathematics},
which introduces class~68V (\textit{Computer science support for mathematical research and practice}) and in particular 68V20 (\textit{Formalization of mathematics}).

Kunen's interest in the area of computational logic dates back to the 1980s.
He published papers on the theory of logic programming \cite{kunen-negation,kunen-signed}
and on resolution theorem proving~\cite{hart-single,kunen-semantics-answer}.
He later became interested in the Boyer--Moore theorem prover, a distinctive semi-automatic system based on a quirky, quantifier free first-order logic; his aim seems to have been to examine the strength of that logic \cite{kunen-ramsey,kunen-nonconstructive}. 
He was also aware of my own work on formalising G\"odel's constructible universe using Isabelle/ZF \cite{paulson-consistency}. 
I'd like to think that he would take an interest in a formalisation of work by \erd{} and Milner~\cite{erdos-theorem-partition} on the partition calculus.

\section{Ordinal partition relations}\label{sec:ordinal_partitions}

\erd{} and Rado introduced the \textit{partition calculus} in 1952 to investigate a family of problems related to Ramsey's theorem~\cite{erdos-partition,Hajnal-Larson}. 
Let $i$, $j$, \ldots{} denote natural numbers while $\alpha$, $\beta$, \ldots{} denote set-theoretic ordinals. Let $[A]^n$ denote the set of $n$-element subsets of a given set~$A$. Write $\tp A$ for the order type of~$A$.

Now we can define \textit{partition notation}: $\alpha\longrightarrow (\beta_0, \ldots,\beta_{k-1})^n$
means that for every partition of the set $[\alpha]^n$  into $k$ parts or ``colours'' $C_0$, \ldots, $C_{k-1}$, there exists $i<k$ and a subset $B\subseteq\alpha$ of order type $\beta_i$ such that $[B]^n\subseteq C_i$. Such a $B$ is said to be
\textit{$i$-monochromatic}. The same notation can be used with order types replaced by cardinalities.

Below we consider only the special case $\alpha\longrightarrow (\beta, \gamma)^2$ and omit the superscript. The negation of $\alpha\longrightarrow (\beta, \gamma)$ is written $\alpha\narrows (\beta, \gamma)$.
In this notation, the infinite Ramsey theorem becomes $\omega\longrightarrow (\omega, \omega)$. 
A straightforward construction proves $\alpha\narrows (|\alpha|+1,\, \omega)$ for $\alpha>\omega$, while
$\alpha\longrightarrow (\alpha,2)$ is trivial. 
If $\alpha$ is not a power of $\omega$ (or zero) then there exist ordinals $\beta$, $\gamma<\alpha$ such that $\alpha=\beta+\gamma$ \cite[p.\ts43]{kunen80},
from which it easily follows that $\alpha\narrows(\alpha,3)$.
These and other facts raise the question~\cite[\S3.2]{erdos-unsolved}, 
for which $m$ and countable ordinals $\alpha$ do we have $\alpha\longrightarrow (\alpha,m)$?

Kunen's  interest in partition theorems is clear in a result he announced in 1971, which is equivalent to $\kappa \longrightarrow (\kappa, \alpha)$:
\begin{quote}
Let $\kappa$ be a real-valued measurable cardinal and $\mu$ a normal measure on~$\kappa$. Let $A\subseteq [\kappa]^2$. Then either (i) there is a subset, $X\subseteq\kappa$, such that $\mu(X)>0$ and $[X]^2\subseteq A$, or (ii) for all countable ordinals~$\alpha$, there is an $X\subseteq\kappa$ such that $X$ has order type~$\alpha$ and $[X]^2\cap A=\emptyset$. The proof uses a generalization of the zero-one law.~\cite{kunen-partition-theorem}
\end{quote}

\section{Introduction to Isabelle}\label{sec:isabelle}

Isabelle is an interactive theorem prover based on a logical framework: a minimal formalism intended for representing formal proofs in a variety of logics~\cite{paulson-found}.
Isabelle/ZF supports first-order logic and set theory, and has been used to formalise the constructible universe~\cite{paulson-consistency} and forcing~\cite{gunther-forcing}.
But its most popular instance by far is Isabelle/HOL \cite{isa-tutorial}, supporting higher-order logic. All versions of Isabelle share a substantial code base, including a sophisticated interactive environment and tools for automatic simplification and logical reasoning.
However, Isabelle/HOL extends all that with specialised, powerful automation for proving theorems and detecting counterexamples~\cite{paulson-from-lcf}.

Isabelle's higher-order logic is an extension of Church's simple type theory~\cite{church40}. It assumes the axiom of choice.
It has basic types such as \isa{nat} (the natural numbers) and \isa{bool} (the truth values, and hence the type of formulas).
It has function types such as \isa{$\alpha$\isasymRightarrow$\beta$} and (postfix) type operators such as {$\alpha$\,set}, sets over type~$\alpha$.
Thus, types can take other types as parameters, but they can't take other values as parameters: there are no ``dependent types''. My colleagues and I are pursuing the thesis that simple type theory is not merely sufficient to formalise mathematics~\cite{bordg-simple-tr} but superior to strong type theories that make automation difficult and introduce complications such as intensional equality.

The set theoretic developments reported here were actually undertaken using Isabelle/HOL, augmented with the ZF axioms; they would have been harder in the more basic proof environment of Isabelle/ZF\@.
The axiomatisation of ZF in HOL~\cite{ZFC_in_HOL-AFP} introduces a type \isa{V}, the type of all ZF sets.
Type \isa{V set} is effectively the type of classes, and any \isa{small} class can be mapped to the corresponding element of~\isa{V}\@. Transfinite recursion is easily obtained from Isabelle/HOL's support of recursive function definitions, and the formal development of set theory includes ordinals, cardinals, order types of well-founded relations, Cantor normal form and other essential material, the proofs mostly taken from Kunen's well-known textbook~\cite{kunen80}. 

An order type is always an ordinal in this formalisation of ZF\@.
In the general case, \isa{ordertype} applies to any set~\isa{A} and well-founded relation~\isa{r}, but here that relation is always set membership. \erd{} and Milner \cite{erdos-theorem-partition} actually considered order types of arbitrary orderings, but the special case of ordinals is sufficient for our application of formalising Larson~\cite{larson-short-proof}.

\section{Outline of the proof}

\erd{} and Milner \cite{erdos-theorem-partition} proved that if $\nu$ is a countable ordinal and $n<\omega$ then 
\begin{equation}
\omega^{1+\nu n}\longrightarrow(2^n, \omega^{1+\nu}). \label{eqn:thm}
\end{equation}
They claim to have known the result since 1959, from Milner's PhD work. Remarkably, the published proof of the main lemma contained so many errors that their five page paper needed a full-page correction~\cite{erdos-theorem-partition-corr}, replacing the core of the original proof. These errors somehow escaped the notice of the authors, the PhD examiners and the original referee. That so many pairs of eyes could overlook so many errors is evidence of the need for more formal scrutiny of published mathematics.

The proof is highly technical, and for the full details, readers should consult the \erd--Milner paper itself~\cite{erdos-theorem-partition} and crucially, the corrections~\cite{erdos-theorem-partition-corr}.
Below we shall simply examine the milestones of the proof, with comments on the special difficulties occasioned by their formalisation.

\erd{} and Milner rely on a more general theorem: that if $\alpha\longrightarrow(k, \gamma)$ and $k\ge2$, then
\[ \alpha\beta\longrightarrow(2k,\,\gamma\vee \omega\beta). \]
Already some complications are evident. In this statement of the theorem, the symbol~$\vee$  extends the partition notation to allow a choice between a 1-monochromatic set of order type $\gamma$ or one of type $\omega\beta$, and the Greek letters refer to order types in the general sense: where two orderings have the same order type when there exists an order-preserving bijection between them. It is not even clear how to formalise this general statement in ZFC, where an order type is a proper class. So the first step is to reformulate the theorem (including the technical condition that $\beta$ is a ``strong type'', and others) for  ordinals: this special case suffices for the main result.

Here is the statement above for the case when $\alpha$, $\beta$ and $\gamma$ range over ordinals, 
$\alpha$ is indecomposable and $\beta$ is countable. 
\begin{equation}
\text{If } \alpha\longrightarrow(k, \gamma)
\text{ and } k\ge2
\text{ then }
 \alpha\beta\longrightarrow(2k,\,\min(\gamma, \omega\beta)). \label{eqn:thm2}
\end{equation}
Because the ordinals are linearly ordered, the choice between finding a set of type~$\gamma$ or a set of type~$\omega\beta$ no longer requires a disjunction but just taking their minimum. 
This statement~(\ref{eqn:thm2}) suffices to prove the original claim, $\omega^{1+\nu n}\longrightarrow(2^n, \omega^{1+\nu})$. \erd{} and Milner provide a full inductive argument, reproduced below with trivial substitutions:
\begin{quote}
Suppose (\ref{eqn:thm}) holds for some integer $n\ge1$. Applying the above theorem with $k=2^n$, $\alpha = \omega^{1+\nu n}$, $\beta=\omega^\nu$, $\gamma=\omega^{1+\nu}$, we see that (\ref{eqn:thm}) also holds with $n$ replaced by $n+1$. Since (\ref{eqn:thm}) holds trivially for $n=0$, it follows that (\ref{eqn:thm}) holds for all $n<\omega$. \cite[p.\ts503]{erdos-theorem-partition}
\end{quote}
This proof was easy to formalise and is presented in full (Fig.\ts\ref{fig:main}). The assumption $n\ge1$ turns out to be unnecessary. Some notes on the syntax: to formalise $\omega^{1+\nu n}\longrightarrow(2^n, \omega^{1+\nu})^2$ requires the constant for the partition relation, \isa{partn\_lst\_VWF}, and an explicit conversion from natural numbers to the corresponding finite ordinals, \isa{ord\_of\_nat}. Key claims are labelled with \isakeyword{shows} and intermediate claims with \isakeyword{have}. Justifications begin with \isakeyword{using} followed by the names of prior results or with \isakeyword{by}, followed by a proof method such as \isa{auto}.
\begin{isabelle}
\ \ \ \ "partn\_lst\_VWF\ (\isasymomega \isasymup (1\ +\ \isasymnu *n))\ [ord\_of\_nat\ (2\isacharcircum n),\ \isasymomega \isasymup (1+\isasymnu )]\ 2"
\end{isabelle}

In the sequel, we shall only be concerned with proving the theorem~(\ref{eqn:thm2}).

\begin{figure*}
\begin{isabelle}
\isacommand{theorem}\ Erdos\_Milner:\isanewline
\ \ \isakeyword{assumes}\ \isasymnu :\ "\isasymnu \ \isasymin \ elts\ \isasymomega 1"\isanewline
\ \ \isakeyword{shows}\ "partn\_lst\_VWF\ (\isasymomega \isasymup (1\ +\ \isasymnu *n))\ [ord\_of\_nat\ (2\isacharcircum n),\ \isasymomega \isasymup (1+\isasymnu )]\ 2"\isanewline
\isacommand{proof}\ (induction\ n)\isanewline
\ \ \isacommand{case}\ 0\isanewline
\ \ \isacommand{then}\ \isacommand{show}\ ?case\isanewline
\ \ \ \ \isacommand{using}\ partn\_lst\_VWF\_degenerate\ [of\ 1\ 2]\ \isacommand{by}\ simp\isanewline
\isacommand{next}\isanewline
\ \ \isacommand{case}\ (Suc\ n)\isanewline
\ \ \isacommand{have}\ "Ord\ \isasymnu "\isanewline
\ \ \ \ \isacommand{using}\ Ord\_\isasymomega 1\ Ord\_in\_Ord\ assms\ \isacommand{by}\ blast\isanewline
\ \ \isacommand{have}\ "1+\isasymnu \ \isasymle \ \isasymnu +1"\isanewline
\ \ \ \ \isacommand{by}\ (simp\ add:\ \isacartoucheopen Ord\ \isasymnu \isacartoucheclose \ one\_V\_def\ plus\_Ord\_le)\isanewline
\ \ \isacommand{then}\ \isacommand{have}\ [simp]:\ "min\ (\isasymomega \ \isasymup \ (1\ +\ \isasymnu ))\ (\isasymomega \ *\ \isasymomega \ \isasymup \ \isasymnu )\ =\ \isasymomega \ \isasymup \ (1+\isasymnu )"\isanewline
\ \ \ \ \isacommand{by}\ (simp\ add:\ \isacartoucheopen Ord\ \isasymnu \isacartoucheclose \ oexp\_add\ min\_def)\isanewline
\ \ \isacommand{have}\ ind:\ "indecomposable\ (\isasymomega \ \isasymup \ (1\ +\ \isasymnu \ *\ ord\_of\_nat\ n))"\isanewline
\ \ \ \ \isacommand{by}\ (simp\ add:\ \isacartoucheopen Ord\ \isasymnu \isacartoucheclose \ indecomposable\_\isasymomega \_power)\isanewline
\ \ \isacommand{show}\ ?case\isanewline
\ \ \isacommand{proof}\ (cases\ "n\ =\ 0")\isanewline
\ \ \ \ \isacommand{case}\ True\isanewline
\ \ \ \ \isacommand{then}\ \isacommand{show}\ ?thesis\isanewline
\ \ \ \ \ \ \isacommand{using}\ partn\_lst\_VWF\_\isasymomega \_2\ \isacartoucheopen Ord\ \isasymnu \isacartoucheclose \ one\_V\_def\ \isacommand{by}\ auto\isanewline
\ \ \isacommand{next}\isanewline
\ \ \ \ \isacommand{case}\ False\isanewline
\ \ \ \ \isacommand{then}\ \isacommand{have}\ "Suc\ 0\ <\ 2\isacharcircum n"\isanewline
\ \ \ \ \ \ \isacommand{using}\ less\_2\_cases\ not\_less\_eq\ \isacommand{by}\ fastforce\isanewline
\ \ \ \ \isacommand{then}\ \isacommand{have}\ "partn\_lst\_VWF\ (\isasymomega \ \isasymup \ (1\ +\ \isasymnu *n)\ *\ \isasymomega \isasymup \isasymnu )\isanewline
\ \ \ \ \ \ \ \ \ \ \ \ \ \ \ \ \ \ [ord\_of\_nat\ (2\ *\ 2\isacharcircum n),\ \isasymomega \ \ \isasymup \ (1+\isasymnu )]\ \ 2"\isanewline
\ \ \ \ \ \ \isacommand{using}\ Erdos\_Milner\_aux\ [OF\ Suc\ ind,\ \isakeyword{where}\ \isasymbeta \ =\ "\isasymomega \isasymup \isasymnu "]\ \isacartoucheopen Ord\ \isasymnu \isacartoucheclose \ \isasymnu \isanewline
\ \ \ \ \ \ \isacommand{by}\ (auto\ simp:\ countable\_oexp)\isanewline
\ \ \ \ \isacommand{then}\ \isacommand{show}\ ?thesis\isanewline
\ \ \ \ \ \ \isacommand{using}\ \isacartoucheopen Ord\ \isasymnu \isacartoucheclose \ \isacommand{by}\ (simp\ add:\ mult\_succ\ mult.assoc\ oexp\_add)\isanewline
\ \ \isacommand{qed}\isanewline
\isacommand{qed}	
\end{isabelle}
\caption{Isabelle/HOL formal proof of the main inductive argument} \label{fig:main}
\end{figure*}

\subsection{Preliminaries}

First, some notation and conventions.
The paper~\cite[p.\ts503]{erdos-theorem-partition} refers to fixed sets $S$ of type $\alpha\beta$ and $B$ of type~$\beta$. But since ordinals are sets, working with ordinals rather than order types allows us to use $\alpha\beta$ and~$\beta$ as the required sets. 

$A<B$ means if $x\in A$ and $y\in B$ then $x<y$.

The theorem (\ref{eqn:thm2}) is trivial unless $\alpha>1$ and $\beta\not=0$, which we assume below. By convention, $A$, $A'$, $A_1$, etc.\ denote subsets of $\alpha\beta$ having type $\alpha$.

\subsection{Every ordinal is strong}
\label{sec:strong}
The property that $\beta$ is a strong type does not have to be assumed because every ordinal is strong, meaning if $D \subseteq \beta$ then there are sets $D_1$, \ldots, $D_n \subseteq D$ such that
\begin{itemize}
	\item $\tp D_i$ is indecomposable for $i = 1$, \ldots,~$n$
	\item if $M \subseteq D$ and $\tp (M \cap D_i) \ge \tp D_i$ for $i = 1$, \ldots,~$n$, then $\tp M = \tp D$.
\end{itemize}
In Isabelle/HOL, the theorem statement looks like this, where \isa{L} is a list of sets and \isa{List.set\ L} stands for the set $\{D_1,\ldots,D_n\}$:
\begin{isabelle}
\isacommand{proposition}\ strong\_ordertype\_eq:\isanewline
\ \ \isakeyword{assumes}\ "D\ \isasymsubseteq \ elts\ \isasymbeta "\ \isakeyword{and}\ "Ord\ \isasymbeta "\isanewline
\ \ \isakeyword{obtains}\ L\ \isakeyword{where}\ "\isasymUnion (List.set\ L)\ =\ D"\isanewline
\ \ \ \ \ \ \ \ \ \ \ \ \ \ \ \ \ \ \ "\isasymAnd X.\ X\ \isasymin \ List.set\ L\ \isasymLongrightarrow \ indecomposable\ (tp\ X)"\isanewline
\ \ \ \ \isakeyword{and}\ "\isasymAnd M.\ \isasymlbrakk M\ \isasymsubseteq \ D;\ \isasymAnd X.\ X\ \isasymin \ List.set\ L\ \isasymLongrightarrow \ tp\ (M\ \isasyminter \ X)\ \isasymge \ tp\ X\isasymrbrakk\isanewline
\ \ \ \ \ \ \ \ \ \ \ \ \ \  \ \isasymLongrightarrow \ tp\ M\ =\ tp\ D"
\end{isabelle}
The proof involves writing $\tp D$ in Cantor normal form:
\[ \tp D = \omega^{\beta_1}\cdot l_1 + \cdots + \omega^{\beta_n}\cdot l_n \]
where $\tp D\ge \beta_1>\cdots>\beta_n$ and $1 \le l_i < \omega$ for $i = 1$, \ldots,~$n$. Through the bijection between $\tp D$ and $D$, this divides $D$ into the desired $D_1$, \ldots, $D_n$.
The proof is straightforward in principle, but somehow the formalisation is 200 lines long.
The paper mentions Cantor normal form \cite[p.\ts502]{erdos-theorem-partition} but gives no other hints.

\subsection{A remark about indecomposable ordinals} \label{sec:remark}

The proof relies on the following observation~\cite{erdos-theorem-partition-corr}:
if $x\in A$ and $A_1\subseteq A$, then there is $A_2\subseteq A_1$ such that $\{x\}<A_2$.

Recalling that $\tp A=\alpha$, consider the bijection~$\phi$ between $A$ and~$\alpha$. Then $\phi(x)<\alpha$ and we can define $A_x=\{y\in A:\phi(y)\le \phi(x)\}$ and $A_2 = A_1 \setminus A_x$. Then $\{x\}<A_2$ by construction and $\tp A_2=\alpha$ follows because $\alpha$ is indecomposable. The formalisation is a fairly routine 60 lines, not difficult but tiresome for a straightforward claim stated without proof. Here is the formal version of the theorem statement:
\begin{isabelle}
\isacommand{proposition}\ indecomposable\_imp\_Ex\_less\_sets:\isanewline
\ \ \isakeyword{assumes}\ "indecomposable\ \isasymalpha "\ \isakeyword{and}\ "\isasymalpha \ >\ 1"\ \isanewline
\ \ \ \ \isakeyword{and}\ "tp\ A\ =\ \isasymalpha "\ "small\ A"\ "A\ \isasymsubseteq \ ON"\isanewline
\ \ \ \ \isakeyword{and}\ "x\ \isasymin \ A"\ \isakeyword{and}\ "tp\ A1\ =\ \isasymalpha "\ "A1\ \isasymsubseteq \ A"\isanewline
\ \ \isakeyword{obtains}\ A2\ \isakeyword{where}\ "tp\ A2\ =\ \isasymalpha "\ "A2\ \isasymsubseteq \ A1"\ "\{x\}\ \isasymlless \ A2"
\end{isabelle}
Note that the keyword \isakeyword{obtains} is a way of expressing an existential conclusion, and that the implicit order types of~$A$, $A_1$, $A_2$ need to be written out.

\section{Proving the theorem}

Recall that our task is to prove that if $\alpha\longrightarrow(k, \gamma)$ and $k\ge2$ then
\[ \alpha\beta\longrightarrow(2k,\,\min(\gamma, \omega\beta)) \]
for ordinals $\alpha$, $\beta$, $\gamma$ where $\alpha$ is indecomposable and $\beta$ is countable. Here is the formal version of the statement above, where \isa{\isasymbeta \ \isasymin \ elts\ \isasymomega 1} means $\beta<\omega_1$.
\begin{isabelle}
\isacommand{theorem}\ Erdos\_Milner\_aux:\isanewline
\ \ \isakeyword{assumes}\ "partn\_lst\_VWF\ \isasymalpha \ [k,\ \isasymgamma ]\ 2"\isanewline
\ \ \ \ \isakeyword{and}\ "indecomposable\ \isasymalpha "\ \isakeyword{and}\ "k\ >\ 1"\ "Ord\ \isasymgamma "\ \isakeyword{and}\ \isasymbeta :\ "\isasymbeta \ \isasymin \ elts\ \isasymomega 1"\isanewline
\ \ \isakeyword{shows}\ "partn\_lst\_VWF\ (\isasymalpha *\isasymbeta )\ [ord\_of\_nat\ (2*k),\ min\ \isasymgamma \ (\isasymomega *\isasymbeta )]\ 2"	
\end{isabelle}

The proof considers the set $[\alpha\beta]^2$ partitioned into sets $K_0$ and $K_1$ by a colouring function $f:[\alpha\beta]^2\to \{0,1\}$. Then either (again paraphrasing the authors~\cite[p.\ts503]{erdos-theorem-partition})
\begin{enumerate}[(i)]
	\item there is $X\in [\alpha\beta]^{2k}$ such that $[X]^2 \subseteq K_0$, or
	\item there is $C\subseteq \alpha\beta$ such that $\tp C = \gamma$ and $[C]^2 \subseteq K_1$, or
	\item there is $Z\subseteq \alpha\beta$ such that $\tp Z = \omega\beta$ and $[Z]^2 \subseteq K_1$.
\end{enumerate}
The proof assumes that (i) and (ii) above are both false and deduces~(iii). Let $\{\gamma_m:m<\omega\}$ be an enumeration of~$\beta$ that repeats every element infinitely often.
The 1-monochromatic set $Z$ is constructed by an elaborate enumeration along with increasing families of sets~$\{A^{(n)}_\nu\}_{\nu<\beta}$ satisfying $\tp(Z\cap A^{(m)}_{\gamma_m})=\omega$ for $m<\omega$, from which it can be shown that $\tp Z=\omega\beta$. The formalisation of the full argument takes nearly a thousand lines, and here we look at some milestones. 

Near the start of the proof, we find the claim
\begin{quote}
(8)\quad If $A\subseteq \alpha\beta$, then there is $X\in[A]^k$ such that $[X]^2\subseteq K_0$. This follows from the hypothesis $\alpha\longrightarrow(k, \gamma)$ and the assumed falsity of~(ii).~\cite[p.\ts503]{erdos-theorem-partition}
\end{quote}
The claim seems obvious enough and no further details are given; its formal proof of about 50 lines could possibly be streamlined with the help of higher-level lemmas about partition relations. The claim is embedded in the main proof:
\begin{isabelle}
\ \ \ \ \isacommand{have}\ Ak0:\ "\isasymexists X\ \isasymin \ [A]\isactrlbsup k\isactrlesup .\ f\ `\ [X]\isactrlbsup 2\isactrlesup \ \isasymsubseteq \ \{0\}"\isanewline
\ \ \ \ \ \ \isakeyword{if}\ "A\ \isasymsubseteq \ elts\ (\isasymalpha *\isasymbeta )"\ \isakeyword{and}\ "tp\ A\ \isasymge \ \isasymalpha "\ \isakeyword{for}\ A	
\end{isabelle}
Here, \isa{f\ `\ [X]\isactrlbsup 2\isactrlesup \ \isasymsubseteq \ \{0\}} expresses $[X]^2\subseteq K_0$ in terms of an image involving the colouring function, $f$.

The next stage of the proof requires a new definition: 
\[ K_i(x)=\{y\in\alpha\beta:\{x,y\}\in K_i\} \]
for $x\in\alpha\beta$ and $i<2$, and the claim is~\cite[p.\ts503]{erdos-theorem-partition}
\begin{quote}
(9)\quad	 Suppose $D\subseteq\beta$, $A_\nu\subseteq\alpha\beta$ (for $\nu\in D$), $A\subseteq\alpha\beta$. For $x\in A$ let
\[ M(x)=\{v\in D:\tp(K_1(x)\cap A_\nu)\ge\alpha\}. \]
Then $\tp\{x\in A:\tp M(x) \ge \tp D\} \ge \alpha$. 
\end{quote}
This claim is first proved in a special (weaker) form, assuming that $\tp D$ is indecomposable, and then in a general form, dropping that assumption. The authors prove the specialised version in half a dozen lines using their claim~(8)  and take a further five lines, using the property that $\beta$ is strong, to achieve the general version. The formal proof of the special version is 70 lines long, including a lengthy calculation, while that for the general form is 40 lines, with an induction on the decomposition of~$\tp D$. Here is the statement of the special version:
\begin{isabelle}
\ \ \ \ \isacommand{have}\ 9:\ "tp\ \{x\ \isasymin \ A.\ tp\ (M\ D\ \isasymAA \ x)\ \isasymge \ tp\ D\}\ \isasymge \ \isasymalpha "\isanewline
\ \ \ \ \ \ \isakeyword{if}\ "indecomposable\ (tp\ D)"\ \isakeyword{and}\ "D\ \isasymsubseteq \ elts\ \isasymbeta "\isanewline
\ \ \ \ \ \ \ \ \isakeyword{and}\ "A\ \isasymsubseteq \ elts\ (\isasymalpha *\isasymbeta )"\ \isakeyword{and}\ "tp\ A\ =\ \isasymalpha "\isanewline
\ \ \ \ \ \ \ \ \isakeyword{and}\ "\isasymAA \ \isasymin \ D\ \isasymrightarrow \ \{X.\ X\ \isasymsubseteq \ elts\ (\isasymalpha *\isasymbeta )\ \isasymand \ tp\ X\ =\ \isasymalpha \}"\ \ \isakeyword{for}\ D\ A\ \isasymAA 	
\end{isabelle}
They continue~\cite[p.\ts504]{erdos-theorem-partition} with a not-quite-trivial instance:
\begin{quote}
As a special case of (9) (with $\tp D=1$) we have:\\
(9')\quad If $A$, $A'\subseteq\alpha\beta$. then $\tp\{x\in A': \tp(K_1(x)\cap A)\ge\alpha\}  \ge \alpha$.	
\end{quote}
The preliminaries of the proof conclude with a claim that follows with the help of (9) and (9') in seven lines of text. The formalisation was tough, 300 lines, but identified some small errors.
\begin{quote}
(10)\quad Let $F$ be a finite subset of~$\beta$ and let $\{A_\nu\}_{\nu<\beta}$, $A$ be such that
$\bigcup_{\nu<\beta} A_\nu\subseteq \alpha\beta$ and $A\subseteq \alpha\beta$. Then there are
$x_0\in A$ and a strictly increasing map $g:\beta\to\beta$ such that $g(\nu)=\nu$ ($\nu\in F$) and $\tp (K_1((x_0)\cap A_{g(\nu)}) \ge \alpha$ ($\nu\in\beta$).
\end{quote}
Here is the corresponding formal statement, where \isa{\isasymAA} is the family~$\{A_\nu\}_{\nu<\beta}$:
\begin{isabelle}
\ \ \ \ \isacommand{have}\ 10:\ "\isasymexists x0\isasymin A.\ \isasymexists g\ \isasymin \ elts\ \isasymbeta \ \isasymrightarrow \ elts\ \isasymbeta .\ strict\_mono\_on\ g\ (elts\ \isasymbeta )\isanewline
\ \ \ \ \ \ \ \ \ \ \ \ \ \ \ \ \ \ \ \ \ \ \ \ \ \ \isasymand\ (\isasymforall \isasymnu \ \isasymin \ F.\ g\ \isasymnu \ =\ \isasymnu )\isanewline
\ \ \ \ \ \ \ \ \ \ \ \ \ \ \ \ \ \ \ \ \ \ \ \ \ \ \isasymand \ (\isasymforall \isasymnu \ \isasymin \ elts\ \isasymbeta .\ tp\ (K\ 1\ x0\ \isasyminter \ \isasymAA \ (g\ \isasymnu ))\ \isasymge \ \isasymalpha )"\isanewline
\ \ \ \ \ \ \isakeyword{if}\ F:\ "finite\ F"\ "F\ \isasymsubseteq \ elts\ \isasymbeta "\isanewline
\ \ \ \ \ \ \ \ \isakeyword{and}\ A:\ "A\ \isasymsubseteq \ elts\ (\isasymalpha *\isasymbeta )"\ "tp\ A\ =\ \isasymalpha "\isanewline
\ \ \ \ \ \ \ \ \isakeyword{and}\ \isasymAA :\ "\isasymAA \ \isasymin \ elts\ \isasymbeta \ \isasymrightarrow \ \{X.\ X\ \isasymsubseteq \ elts\ (\isasymalpha *\isasymbeta)\ \isasymand \ tp\ X\ =\ \isasymalpha \}"\isanewline
\ \ \ \ \ \ \isakeyword{for}\ F\ A\ \isasymAA 
\end{isabelle}
One source of difficulty is that the proof expresses $\beta$ as a union of increasing sets $D_0<\{\nu_1\}<D_1<\cdots<\{\nu_p\}<D_p$, where $F=\{\nu_1,\ldots,\nu_p\}$:
\begin{equation}
\beta = D_0\cup\{\nu_1\}\cup D_1\cup\cdots\cup \{\nu_p\}\cup D_p, \label{eqn:beta}	
\end{equation}
Our intuition that $F$ ``obviously'' cuts $\beta$ into segments does not help; every detail of the relationship between $\beta$ and the $D_i$ must be formalised. 

Now the stage is set for the construction of the required 1-monochromatic set $Z=\{x_0,x_1,\ldots,x_n,\ldots\}$ and the increasing family $A^{(n)}_\nu$, by complete induction on~$n$. Here we must switch to the corrigendum~\cite{erdos-theorem-partition-corr}, which replaces claims (12)--(17) of the original proof. The ``remark'' noted above (\S\ref{sec:remark}) and claim~(10) are used here to yield $x_n\in A$, a strictly increasing map $g_n:\beta\to\beta$ and sets $\{A^{(n+1)}_\nu\}_{\nu<\beta}$ satisfying certain properties. The formalisation of this section is about 360 lines. These inductive constructions are typical of Ramsey theory, and better ways of formalising them are needed. However, many of the technical complications are inherent in the mathematics itself.

\section{Observations and Conclusions}

The task for formalising \erd{} and Milner's paper arose in the context of a larger project, with D\v zamonja and Koutsoukou-Argyraki~\cite{dzamonja-formalising-arxiv}, to formalise Larson's proof~\cite{larson-short-proof} that $\omega^\omega\longrightarrow(\omega^\omega, m)$ for all $m<\omega$. Her proof relies on 
\begin{equation}
\omega^{n \cdot k} \longrightarrow (\omega^n, k),	\label{eqn:erdos}
\end{equation}
which she credits to \erd{}, remarking that the Erd\H{o}s-Milner paper ``has a sharper result''~\cite[p.\ts134]{larson-short-proof}.
The claim~(\ref{eqn:erdos}) is trivial if $k=0$ or $n=0$, and otherwise put $\nu=n-1$ and $n=k-1$ in
$\omega^{1+\nu n}\longrightarrow(2^n, \omega^{1+\nu})$. This is a routine calculation and the formalisation is just 30 lines long.

The revision control logs hold the detailed history of the formal development of \erd{}--Milner~\cite{erdos-theorem-partition,erdos-theorem-partition-corr}. I started scrutinising the paper on 3 February 2020. By 7 Feburary, I had proved that the Theorem implied the paper's headline result. By 12 February, I had proved (8) and had started~(9).%
\footnote{The equation numbers here refer to \erd{} and Milner~\cite{erdos-theorem-partition}.}
The general case of the latter required first proving the lemma about strong types (\S\ref{sec:strong}), and it was not until 24 February that I managed to prove~(9) and the corollary (9'). By 1 March, I had proved (10).
Next on the agenda was the inductive construction of $Z$ and the $A^{(n)}_\nu$ families, which was a struggle. The logs refer to a number of unsuccessful attempts to do a ``big induction''.
The entry for 12 March says ``replaced the big induction by a primitive recursion setup'' and that (12)--(16) had been proved. That would have used the remark discussed in~\S\ref{sec:remark} above and included the construction of the 1-monochromatic set~Z\@. The last step was to establish $Z$'s order type; (17), (18) and the necessary half%
\footnote{That is, $\tp(Z\cap A^{(m)}_{\gamma_m})\ge\omega$; the other direction wasn't clear to me.}
of (19) were done on successive days. The formalisation was complete by March 17. 
It can be found online~\cite{Ordinal_Partitions-AFP}, within the formalisation of Larson.

Was it worth it?
The process took 44 days, during the middle of a normal University term with the usual schedule of teaching and administrative duties. Recalling that the paper was only five pages long, this equates to nine days per page to understand and formalise the material. And while the headline result was obtained in full generality, it was proved from a theorem about order types which here was formalised only for the special case of ordinals. The effort required to formalise mathematics is surely still prohibitive. Nevertheless, the numerous errors in the original paper are a reminder of how easy it is to make mistakes, especially perhaps in combinatorial proofs. Every formalisation of a technically difficult piece of mathematics provides strong assurance of its correctness. In the case of Isabelle/HOL, it also often yields a document that is readable enough to be examined by anybody who still has doubts.

Every time we formalise something in a new area of mathematics, we discover certain things that particularly hard to formalise. Formalisation will never catch up with mathematical intuition. Many obvious deductions---such as the partition of~$\beta$, our~(\ref{eqn:beta}) above---seem to take a disproportionate amount of work, and our only consolation is that some obvious deductions are false. The inductive constructions typical of Ramsey and partition theory are an area where we need further work to find more compact, natural and readable formal proofs.

\paragraph*{Acknowledgements}
Thanks to Mirna D{\v z}amonja (who proposed the project in the first place) and to Angeliki Koutsoukou-Argyraki for discussions and help. The ERC supported this research through the Advanced Grant ALEXANDRIA (Project GA 742178).  

\section*{References}

\bibliographystyle{plain}

\bibliography{string,atp,general,isabelle,theory,crossref}

\end{document}